# ON SOLVING GENERAL LINEAR EQUATIONS IN THE SET OF NATURAL NUMBERS

Florentin Smarandache
University of New Mexico
200 College Road
Gallup, NM 87301, USA
E-mail: smarand@unm.edu

The utility of this article is that it establishes if the number of the natural solutions of a general linear equation is limited or not. We will show also a method of solving, using integer numbers, the equation $ax - by = c$ (which represents a generalization of lemmas 1 and 2 of [4]), an example of solving a linear equation with 3 unknowns in N, and some considerations on solving, using natural numbers, equations with $n$ unknowns.

Let's consider the equation:

(1) $\quad \sum_{i=1}^{n} a_i x_i = b$ with all $a_i, b \in \mathbb{Z}$, $a_i \neq 0$, and the greatest common factor $(a_1, ..., a_n) = d$.

**Lemma 1**: The equation (1) admits at least a solution in the set of integers, if $d$ divides $b$.

This result is classic.

In (1), one does not diminish the generality by considering $(a_1, ..., a_n) = 1$, because in the case when $d \neq 1$, one divides the equation by this number; if the division is not an integer, then the equation does not admit natural solutions.

It is obvious that each homogeneous linear equation admits solutions in $\mathbb{N}$: at least the banal solution!

## PROPERTIES ON THE NUMBER OF NATURAL SOLUTIONS OF A GENERAL LINEAR EQUATION

We will introduce the following definition:

**Definition 1**: The equation (1) has variations of sign if there are at least two coefficients $a_i, a_j$ with $1 \leq i, j \leq n$, such that $\text{sign}(a_i \cdot a_j) = -1$

**Lemma 2**: An equation (1) which has sign variations admits an infinity of natural solutions (generalization of lemma 1 of [4]).

*Proof*: From the hypothesis of the lemma it results that the equation has $h$ no null positive terms, $1 \leq h \leq n$, and $k = n - h$ non null negative terms. We have $1 \leq k \leq n$; it is supposed that the first $h$ terms are positive and the following $k$ terms are negative (if not, we rearrange the terms).

We can then write:

$$\sum_{t=1}^{h} a_t x_t - \sum_{j=h+1}^{n} a'_j x_j = b \text{ where } a'_j = -a_j > 0.$$

Let's consider $0 < M = [a_1, ..., a_n]$, the least common multiple, and $c_i = |M/a_i|$, $i \in \{1, 2, ..., n\}$.

Let's also consider $0 < P = [h, k]$, the least common multiple, and $h_1 = P/h$ and $k_1 = P/k$.

Taking $\begin{cases} x_t = h_1 c_t \cdot z + x_t^0, & 1 \leq t \leq h \\ x_j = k_1 c_j \cdot z + x_j^0, & h+1 \leq j \leq n \end{cases}$

where $z \in \mathbb{N}$, $z \geq \max\left\{\left[\dfrac{-x_t^0}{h_1 c_t}\right], \left[\dfrac{x_j^0}{k_1 c_j}\right]\right\} + 1$ with $[\gamma]$ meaning integer part of $\gamma$, i.e. the greatest integer less than or equal to $\gamma$, and $x_i^0$, $i \in \{1, 2, ..., n\}$, a particular integer solution (which exists according to lemma 1), we obtain an infinity of solutions in the set of natural numbers for the equation (1).

**Lemma 3**:
  a) An equation (1) which does not have variations of sign has at maximum a limited number of natural solutions.
  b) In this case, for $b \neq 0$, constant, the equation has the maximum number of solutions if and only if all $a_i = 1$ for $i \in \{1, 2, ..., n\}$.

*Proof*: (see also [6]).
  a) One considers all $a_i > 0$ (otherwise, multiply the equation by -1).
If $b < 0$, it is obvious that the equation does not have any solution (in $\mathbb{N}$).
If $b = 0$, the equation admits only the trivial solution.
If $b > 0$, then each unknown $x_i$ takes positive integer values between 0 and $b/a_i = d_i$ (finite), and not necessarily all these values. Thus the maximum number of solutions is lower or equal to: $\prod_{i=1}^{n}(1 + d_i)$, which is finite.

  b) For $b \neq 0$, constant, $\prod_{i=1}^{n}(1 + d_i)$ is maximum if and only if $d_i$ are maximum, i.e. iff $a_i = 1$ for all $i$, where $i = \{1, 2, ..., n\}$.

**Theorem 1**. The equation (1) admits an infinity of natural solutions if and only if it has variations of sign.
This naturally follows from the previous results.

**Method of solving.**

**Theorem 2.** Let's consider the equation with integer coefficients $ax - by = c$, where $a$ and $b > 0$ and $(a, b) = 1$. Then the general solution in natural numbers of this equation is:
$\begin{cases} x = bk + x_0 \\ y = ak + y_0 \end{cases}$ where $(x_0, y_0)$ is a particular integer solution of the equation,

and $k \geq \max\{[-x_0/b], [-y_0/a]\}$ is an integer parameter (generalization of lemma 2 of [4]).

*Proof:* It results from [1] that the general integer solution of the equation is
$$\begin{cases} x = bk + x_0 \\ y = ak + y_0 \end{cases}$$ where $(x_0, y_0)$ is a particular integer solution of the equation and $k \in \mathbb{Z}$. Since $x$ and $y$ are natural integers, it is necessary for us to impose conditions for $k$ such that $x \geq 0$ and $y \geq 0$, from which it results the theorem.

WE CONCLUDE!

To solve **in the set of natural numbers** a linear equation with $n$ unknowns we will use the previous results in the following way:

a) If the equation does not have variations of sign, because it has a limited number of natural solutions, the solving is made by tests (see also [6])

b) If it has variations of sign and if $b$ is divisible by $d$, then it admits an infinity of natural solutions. One finds its general integer solution (see [2], [5]);

$$x_i = \sum_{j=1}^{n-1} \alpha_{ij} k_j + \beta_i, \ 1 \leq i \leq n$$ where all the $\alpha_{ij}, \beta_i \in \mathbb{Z}$ and the $k_j$ are integer parameters.

By applying the restriction $x_i \geq 0$ for $i$ from $\{1, 2, ..., n\}$, one finds the conditions which must be satisfied by the integer parameters $k_j$ for all $j$ of $\{1, 2, ..., n-1\}$. (c)

The case $n = 2$ and $n = 3$ can be done by this method, but when $n$ is bigger, the conditions (c) become more and more difficult to find.

**Example:** Solve in $\mathbb{N}$ the equation $3x - 7y + 2z = -18$.

Solution: In $\mathbb{Z}$ one obtains the general integer solution:
$$\begin{cases} x = k_1 \\ y = k_1 + 2k_2 \\ z = 2k_1 + 7k_2 - 9 \end{cases}$$ with $k_1$ and $k_2$ in $\mathbb{Z}$.

From the conditions (c) result the inequalities $x \geq 0$, $y \geq 0$, $z \geq 0$. It results that $k_1 \geq 0$ and also:

$k_2 \geq [-k_1/2] + 1$ if $-k_1/2 \notin \mathbb{Z}$, or $k_2 \geq -k_1/2$ if $-k_1/2 \in \mathbb{Z}$;

and $k_2 \geq [(9 - 2k_1)/7] + 1$ if $(9-2k_1)/7 \notin \mathbb{Z}$, or $k_2 \geq (9-2k_1)/7$ if $(9-2k_1)/7 \in \mathbb{Z}$;

that is $k_2 \geq [(2 - 2k_1)/7] + 2$ if $(2-2k_1)/7 \notin \mathbb{Z}$, or $k_2 \geq (2-2k_1)/7 + 1$ if $(2-2k_1)/7 \in \mathbb{Z}$.

With these conditions on $k_1$ and $k_2$ we have the general solution in natural numbers of the equation.